\documentclass[twocolumn]{IEEEtran} 
\usepackage{amsmath}
\usepackage{amssymb}
\usepackage{amsfonts}
\usepackage{epsfig,subfigure}
\usepackage{multirow}

\usepackage{fancyhdr} 
\newcommand{\secref}[1]{Section~\ref{#1}}

\newcommand{\Matrix}[1]{\mathbf{#1} {}}
\newcommand{\Column}[1]{\mathbf{#1} {}}

\newcommand{\IMat}{\Matrix{I}}

\newcommand{\vCol}{\Column{v}}

\newcommand{\AMat}{\Matrix{A}}
\newcommand{\BMat}{\Matrix{B}}

\newcommand{\uCol}{\Column{u}}

\newcommand{\ieeeap}{{\it IEEE Trans. Antennas Propagat.}}
\newcommand{\ieeemtt}{{\it IEEE Trans. Microwave Theory Tech.}}
\newcommand{\motl}{{\it Microwave Opt. Technol. Lett.}}

\newcommand{\el}{{\it Electron. Lett.}}

\newcommand{\myvol}[1]{vol.~{#1}}

\newcommand{\ieeemwcl}{{\it IEEE Microw. Wireless Comp. Lett.}}
\newcommand{\ieeeawpl}{{\it IEEE Antennas. Wireless Propagat. Lett.}}
\newcommand{\half}{\frac{1}{2}}
\newcommand{\quar}{\frac{1}{4}}
\newcommand{\tquar}{\frac{3}{4}}
\newcommand{\oquar}{1\quar}

\newcommand{\uTldCol}{\Column{\tilde{u}}}
\newcommand{\fracdt}[2][]{\frac{#1 \Delta t}{#2}}

\newcommand{\ETld}{\tilde{E}}
\newcommand{\HTld}{\tilde{H}}
\newcommand{\delTld}{\tilde{\partial}}
\renewcommand{\delTld}{\partial}

\begin{document}

\title{
Fundamental Schemes for Efficient Unconditionally Stable Implicit Finite-Difference Time-Domain Methods
}

\author{Eng Leong Tan, {\it Senior Member, IEEE} 
\thanks{The author is with the
School of Electrical and Electronic Engineering,
Nanyang Technological University,
Singapore 639798
(e-mail: eeltan@ntu.edu.sg).}
}

\maketitle
\thispagestyle{fancy}

\begin{abstract}
This paper presents the generalized formulations of 
fundamental schemes 
for efficient unconditionally stable 
implicit finite-difference time-domain (FDTD) methods. 
The fundamental schemes constitute 
a family of implicit schemes that feature 
similar fundamental updating structures, 
which are in simplest forms 
with most efficient right-hand sides. 
The formulations of fundamental schemes 
are presented in terms of generalized matrix operator equations 
pertaining to some classical splitting formulae, 
including those of alternating direction implicit, 
locally one-dimensional and split-step schemes. 
To provide further insights into the implications 
and significance of fundamental schemes, 
the analyses are also extended to many other schemes 
with distinctive splitting formulae. 
Detailed algorithms are described 
for new efficient implementations of 
the unconditionally stable implicit FDTD methods 
based on the fundamental schemes. 
A comparative study of various implicit schemes 
in their original and new implementations is carried out,  
which includes comparisons of their computation costs 
and efficiency gains.  
\end{abstract}



\section{Introduction}

The finite-difference time-domain (FDTD) method 
has been widely used to obtain the numerical solutions of 
Maxwell's equations for investigating electromagnetic wave 
radiation, propagation and scattering problems 
\cite{TafloveBook05}.
For the conventional explicit FDTD method \cite{Yee}, 
the computational efficiency is restricted by 
the Courant-Friedrichs-Lewy (CFL) stability condition, 
which imposes a maximum constraint on the time step size
depending on the spatial mesh sizes.  
To remove the CFL condition, 
the unconditionally stable FDTD method based on 
the alternating direction implicit (ADI) technique 
has been developed \cite{ADI00}, \cite{Namiki00}. 
The success of ADI-FDTD method has brought about 
a resurgence of interest in the unconditionally stable schemes  
not only within the electromagnetics community, 
but also among other scientists and mathematicians at large. 
With its unconditional stability advantage, the ADI-FDTD method 
has been extensively analyzed, improved and extended 
for many applications. 
Since the literature in these aspects is so voluminous, 
we could have easily missed out many references 
that are not directly relevant.  
For more comprehensive survey, we refer the readers 
to the bibliography in \cite{TafloveBook05}, 
\cite{timedomainBook}, \cite{higherorderBook} 
(see also those cited in the author's previous works). 

The idea of ADI scheme  
that has been adapted in the recent celebrated ADI-FDTD, 
can be traced back to the early classic works 
by Peaceman and Rachford \cite{ADI55}. 
Apart from the ADI scheme, many alternative implicit schemes
have also been introduced by researchers in applied mathematics 
to deal with various parabolic, elliptic and hyperbolic 
partial differential systems.  
Some of such schemes that are more common 
and closely related to the present scope 
include locally one-dimensional (LOD) 
and Crank-Nicolson schemes 
\cite{LOD62}-\cite{ThomasBook}.
These schemes have been adapted as well for FDTD solutions of 
Maxwell's equations, leading to unconditionally stable 
LOD-FDTD method \cite{LOD05}-\cite{LOD3D}, 
split-step FDTD approach \cite{splitstep03}-\cite{FuSS04}  
and other Crank-Nicolson-based approximation methods \cite{CN06}. 
Most of these unconditionally stable FDTD methods  
have been built upon the classical implicit schemes 
by adopting their respective splitting formulae directly.  

Despite the successful adaptation of various classical schemes, 
continuing efforts are underway to devise new stable methods 
that are more efficient and simpler to implement. 
Recently, a new efficient algorithm has been presented 
for the ADI-FDTD method \cite{ADIefficient}. 
The algorithm involves updating equations whose right-hand sides 
are much simpler and more concise than 
the conventional implementation \cite{ADI00}, \cite{Namiki00}.  
This leads to substantial reduction 
in the number of arithmetic operations 
required for their computations. 
While the underlying principle of the new algorithm 
has helped make the ADI-FDTD simpler and more efficient, 
it actually has greater significance in its own right. 
This paper aims to clarify and extend such principle 
to arrive at a series of new efficient algorithms 
for various other unconditionally stable implicit FDTD methods, 
including some of those mentioned above. 
Moreover, it will be found that the resultant algorithms 
constitute a family of implicit schemes, 
all of which feature similar fundamental updating structures 
that are in simplest forms with most efficient right-hand sides.  

The organization of this paper is as follows. 
\secref{formulation} presents the formulations of 
fundamental schemes 
in terms of generalized matrix operator equations 
pertaining to some classical splitting formulae, 
including those of ADI, LOD and split-step schemes. 
To provide further insights into the implications 
and significance of fundamental schemes, 
the analyses are also extended to many other schemes 
with distinctive splitting formulae. 
In \secref{implementation}, detailed algorithms 
are described based on the fundamental schemes 
for new efficient implementations of 
the unconditionally stable implicit FDTD methods, 
e.g. ADI-FDTD and LOD-FDTD. 
In \secref{comparison}, 
a comparative study of various implicit schemes 
in their original and new implementations is carried out,  
which includes comparisons of their computation costs 
and efficiency gains.  
The fundamental nature of new implementations 
will also become evident through the comparisons and discussions.


\section{Generalized Formulations}
\label{formulation}

In this section, we present the generalized formulations of 
fundamental schemes for implicit finite-difference methods. 
Starting from some classical implicit schemes, 
their generalized matrix operator equations are revisited and 
reformulated in the simplest and most efficient forms. 
These new forms feature convenient matrix-operator-free right-hand sides 
with least number of terms, which will lead to 
coding simplification in their algorithm implementations. 
For simplicity and clarity, we omit the nonhomogeneous terms 
that appear merely as vectors (without involving matrices).

\subsection{ADI}

The ADI scheme, originated by Peaceman and Rachford \cite{ADI55}, 
is one of the most popular implicit finite-difference schemes 
in use today. 
This scheme calls for generalized splitting formulae in the form 
\begin{subequations}
\label{clsADI}
\begin{gather}
   \Big( \IMat - \fracdt{2} \AMat \Big) \uCol^{n+\half} 
   = \Big( \IMat + \fracdt{2} \BMat \Big) \uCol^{n} 
   \label{clsADIa}
\\
   \Big( \IMat - \fracdt{2} \BMat \Big) \uCol^{n+1}
   = \Big( \IMat + \fracdt{2} \AMat \Big) \uCol^{n+\half} 
   \label{clsADIb}
.\end{gather}
\end{subequations}
For many decades, such splitting formulae form 
the basis of many other numerical methods, 
which include the recent unconditionally stable ADI-FDTD method 
\cite{ADI00, Namiki00}. 
Note that the specific matrix operators $\AMat$ and $\BMat$ of \cite{ADI00, Namiki00} 
(for 3-D Maxwell) are different from those of \cite{ADI55} (for 2-D parabolic), 
even though they conform to the same generalized form \eqref{clsADI}. 
In this section, we shall let these matrices be general 
but would caution that one must choose their operators properly 
for a particular scheme to stay unconditionally stable \cite{CN06}. 
To implement the ADI algorithm, 
it is more convenient to introduce auxiliary variables 
for denoting the right-hand sides of implicit equations. 
This allows us to rewrite the original algorithm as 
\begin{subequations}
\label{orgADI}
\begin{gather}
   \vCol^{n} 
   = \Big( \IMat + \fracdt{2} \BMat \Big) \uCol^{n} 
   \label{orgADIa}
\\
   \Big( \IMat - \fracdt{2} \AMat \Big) \uCol^{n+\half} 
   = \vCol^{n}
   \label{orgADIb}
\\
   \vCol^{n+\half} 
   = \Big( \IMat + \fracdt{2} \AMat \Big) \uCol^{n+\half} 
   \label{orgADIc}
\\      
   \Big( \IMat - \fracdt{2} \BMat \Big) \uCol^{n+1}
   = \vCol^{n+\half}
   \label{orgADId}
\end{gather}
\end{subequations}
where the $\vCol$'s serve as the auxiliary variables. 

If one exploits these auxiliary variables, 
the original algorithm above can be modified into 
one of the more efficient scheme. 
In particular, based on \eqref{orgADId} at one time step backward:  
\begin{gather}
   \vCol^{n-\half}
   = \Big( \IMat - \fracdt{2} \BMat \Big) \uCol^{n}
   \label{proveADIi} 
,\end{gather}
it follows that $\vCol^{n}$ of \eqref{orgADIa} 
is reducible to 
\begin{subequations}
\begin{align}
   \vCol^{n}
   & = \Big( \IMat + \fracdt{2} \BMat \Big) \uCol^{n} 
\\ 
   & = 2 \uCol^{n} - \Big( \IMat - \fracdt{2} \BMat \Big) \uCol^{n} 
\\
   & = 2 \uCol^{n} - \vCol^{n-\half}
   \label{proveADIa}
.\end{align}
\end{subequations}
Furthermore, upon recognizing \eqref{orgADIb}, 
$\vCol^{n+\half}$ of \eqref{orgADIc} 
is also reducible to  
\begin{subequations}
\begin{align}
   \vCol^{n+\half}
   & = \Big( \IMat + \fracdt{2} \AMat \Big) \uCol^{n+\half} 
\\ 
   & = 2 \uCol^{n+\half} - \Big( \IMat - \fracdt{2} \AMat \Big) \uCol^{n+\half} 
\\
   & = 2 \uCol^{n+\half} - \vCol^{n}
   \label{proveADIc}
.\end{align}
\end{subequations}
With \eqref{proveADIa} and \eqref{proveADIc}, 
algorithm \eqref{orgADI} becomes more efficient as 
\begin{subequations}
\label{fastADI}
\begin{gather}
   \vCol^{n} 
   = 2 \uCol^{n} - \vCol^{n-\half} 
   \label{fastADIa}
\\
   \Big( \IMat - \fracdt{2} \AMat \Big) \uCol^{n+\half} 
   = \vCol^{n}
   \label{fastADIb}
\\
   \vCol^{n+\half} 
   = 2 \uCol^{n+\half} - \vCol^{n} 
   \label{fastADIc}
\\      
   \Big( \IMat - \fracdt{2} \BMat \Big) \uCol^{n+1}
   = \vCol^{n+\half}
   \label{fastADId}
.\end{gather}
\end{subequations}

Through a simple re-definition of field variables
\begin{gather}
   \uTldCol^{n} = 2 \uCol^{n}
, \quad
   \uTldCol^{n+\half} = 2 \uCol^{n+\half}
   \label{uTldCol}
,\end{gather}
we may reduce the algorithm further into \cite{ADIefficient} 
\begin{subequations}
\label{fasterADI}
\begin{gather}
   \vCol^{n} 
   = \uTldCol^{n} - \vCol^{n-\half} 
   \label{fasterADIa}
\\
   \Big( \half \IMat - \fracdt{4} \AMat \Big) \uTldCol^{n+\half} 
   = \vCol^{n}
   \label{fasterADIb}
\\
   \vCol^{n+\half} 
   = \uTldCol^{n+\half} - \vCol^{n} 
   \label{fasterADIc}
\\      
   \Big( \half \IMat - \fracdt{4} \BMat \Big) \uTldCol^{n+1}
   = \vCol^{n+\half}
   \label{fasterADId}
.\end{gather}
\end{subequations}
This algorithm proceeds in terms of $\uTldCol$'s for its main iterations, 
and only when the field output is required, one may retrieve the data from 
\begin{gather}
   \uCol^{n+1} = \half \uTldCol^{n+1} 
   \label{fasterADIo}
.\end{gather}
Equations \eqref{fasterADIa}-\eqref{fasterADId} constitute 
the most efficient ADI scheme that has the simplest right-hand sides 
without involving any explicit matrix operator. 
For nonzero initial fields, the algorithm takes 
the input initialization  
\begin{gather}
   \vCol^{-\half} = \Big( \IMat - \fracdt{2} \BMat \Big) \uCol^{0}
   \label{fastADIi}
.\end{gather}
Note that such initialization will not degrade 
the accuracy of ADI scheme since it simply corresponds to 
\eqref{proveADIi} with $n=0$ or \eqref{orgADId} with $n=-1$.  
Furthermore, \eqref{proveADIi} along with \eqref{uTldCol}-\eqref{fasterADIb} 
are just \eqref{orgADIa}-\eqref{orgADIb} or \eqref{clsADIa}. 
Hence the equivalence of both present and original (generalized) 
Peaceman-Rachford schemes becomes evident here.

\subsection{LOD1/SS1}

The LOD scheme, introduced in the early Russian literature 
(cf. \cite{LOD62}, \cite{MitchellBook}), 
is another classical scheme that has been used extensively. 
This scheme calls for generalized splitting formulae in the form 
\begin{subequations}
\label{orgLOD}
\begin{gather}
   \Big( \IMat - \fracdt{2} \AMat \Big) \uCol^{n+\half} 
   = \Big( \IMat + \fracdt{2} \AMat \Big) \uCol^{n} 
   \label{orgLODa}
\\
   \Big( \IMat - \fracdt{2} \BMat \Big) \uCol^{n+1}
   = \Big( \IMat + \fracdt{2} \BMat \Big) \uCol^{n+\half} 
   \label{orgLODb}
.\end{gather}
\end{subequations}
Such splitting formulae have been adopted to develop 
the recent unconditionally stable LOD-FDTD method \cite{LOD05} 
and split-step approach \cite{splitstep03}. 
When $\AMat$ and $\BMat$ do not commute, the scheme is accurate to first order in time.  
Thus, the resultant LOD-FDTD method and split-step approach 
may be denoted by LOD1 and SS1 respectively.  

For most efficiency, the original LOD scheme can be modified  
into simplest form with matrix-operator-free right-hand sides. 
From \eqref{orgLODa}, we have 
\begin{subequations}
\begin{align}
   \Big( \IMat - \fracdt{2} \AMat \Big) \uCol^{n+\half} 
   & = \Big( \IMat + \fracdt{2} \AMat \Big) \uCol^{n} 
\\
   & = 2 \uCol^{n} - \Big( \IMat - \fracdt{2} \AMat \Big) \uCol^{n} 
.\end{align}
\end{subequations}
This can be manipulated readily to give 
\begin{gather}
   \Big( \IMat - \fracdt{2} \AMat \Big) \Big( \uCol^{n+\half} + \uCol^{n} \Big) 
   = 2 \uCol^{n}
\end{gather}
where the vector terms in bracket may be denoted by auxiliary variable 
\begin{gather}
   \vCol^{n+\half} = \uCol^{n+\half} + \uCol^{n} 
.\end{gather}
Similar manipulation applies to \eqref{orgLODb} 
which leads to auxiliary variable $\vCol^{n+1}$. 
Combining all auxiliary and field variables, we arrive at 
\begin{subequations}
\label{fastLOD}
\begin{gather}
   \Big( \half \IMat - \fracdt{4} \AMat \Big) \vCol^{n+\half} 
   = \uCol^{n}
   \label{fastLODa}
\\
   \uCol^{n+\half} 
   = \vCol^{n+\half} - \uCol^{n} 
   \label{fastLODb}
\\
   \Big( \half \IMat - \fracdt{4} \BMat \Big) \vCol^{n+1}
   = \uCol^{n+\half}
   \label{fastLODc}
\\
   \uCol^{n+1} 
   = \vCol^{n+1} - \uCol^{n+\half} 
   \label{fastLODd}
.\end{gather}
\end{subequations}
In these equations, all their right-hand sides are seen 
to be in the most convenient matrix-operator-free form. 
Furthermore, there is no special input initialization required and 
the output field solution is directly available from \eqref{fastLODd}.

\subsection{SS2}

While the former SS1 approach is only first-order accurate in time, 
the general split-step approach actually permits simple extensions 
to achieve higher-order temporal accuracy \cite{splitstep03}. 
Consider the split-step approach of second-order accuracy denoted by SS2.  
It is based on the Strang splitting formulae \cite{Strang68} 
and involves three updating procedures 
\begin{subequations}
\label{orgSS2}
\begin{gather}
   \Big( \IMat - \fracdt{4} \AMat \Big) \uCol^{n+\quar} 
   = \Big( \IMat + \fracdt{4} \AMat \Big) \uCol^{n} 
\\
   \Big( \IMat - \fracdt{2} \BMat \Big) \uCol^{n+\tquar}
   = \Big( \IMat + \fracdt{2} \BMat \Big) \uCol^{n+\quar} 
\\
   \Big( \IMat - \fracdt{4} \AMat \Big) \uCol^{n+1} 
   = \Big( \IMat + \fracdt{4} \AMat \Big) \uCol^{n+\tquar} 
.\end{gather}
\end{subequations}

In line with the previous subsection, 
the algorithm implementation can be made much simpler 
along with improved efficiency by resorting to 
\begin{subequations}
\label{fastSS2}
\begin{gather}
   \Big( \half \IMat - \fracdt{8} \AMat \Big) \vCol^{n+\quar} 
   = \uCol^{n}
   \label{fastSS2a}
\\
   \uCol^{n+\quar} 
   = \vCol^{n+\quar} - \uCol^{n} 
\\
   \Big( \half \IMat - \fracdt{4} \BMat \Big) \vCol^{n+\tquar}
   = \uCol^{n+\quar}
\\
   \uCol^{n+\tquar} 
   = \vCol^{n+\tquar} - \uCol^{n+\quar} 
\\      
   \Big( \half \IMat - \fracdt{8} \AMat \Big) \vCol^{n+1} 
   = \uCol^{n+\tquar}
\\
   \uCol^{n+1} 
   = \vCol^{n+1} - \uCol^{n+\tquar} 
   \label{fastSS2f}
.\end{gather}
\end{subequations}
Again, we have the right-hand sides of these equations 
cast in the most convenient matrix-operator-free form. 
However, since there are three implicit equations to be dealt with, 
the scheme still involves more arithmetic operations 
than those having merely two implicit equations 
with matrix-operator-free right-hand sides.

\subsection{LOD2}

An alternative LOD-FDTD method denoted by LOD2  
has been presented, which does not exhibit the non-commutativity error 
and preserves the second-order temporal accuracy \cite{LOD3D}.
In essence, the main iterations proceed with two updating procedures 
in the same way like LOD1, cf.~\eqref{orgLODa}-\eqref{orgLODb}. 
However, the time indices are associated with advancement of time steps 
from $n+1/4$ to $n+3/4$ and from $n+3/4$ to $n+1(1/4)$ as 
\begin{subequations}
\label{orgLOD2}
\begin{gather}
   \Big( \IMat - \fracdt{2} \AMat \Big) \uCol^{n+\tquar}
   = \Big( \IMat + \fracdt{2} \AMat \Big) \uCol^{n+\quar} 
\\
   \Big( \IMat - \fracdt{2} \BMat \Big) \uCol^{n+\oquar} 
   = \Big( \IMat + \fracdt{2} \BMat \Big) \uCol^{n+\tquar} 
.\end{gather}
\end{subequations}
To relate the fields to those normally at integer time steps, 
such association will prompt for additional processing 
for the input field data 
\begin{gather}
   \Big( \IMat - \fracdt{4} \BMat \Big) \uCol^{\quar} 
   = \Big( \IMat + \fracdt{4} \BMat \Big) \uCol^{0} 
   \label{orgLOD2i}
\end{gather}
as well as for the output field data 
\begin{gather}
   \Big( \IMat + \fracdt{4} \BMat \Big) \uCol^{n+1} 
   = \Big( \IMat - \fracdt{4} \BMat \Big) \uCol^{n+\oquar} 
   \label{orgLOD2o}
.\end{gather} 
Note that the input processing 
is to be invoked only once at the beginning,  
while for the output processing at $n+1$, 
it may be performed separately and independently in parallel 
without disrupting the main iterations. 
Furthermore, one often does not need to frequently output 
all field components, except probably one or two of interest 
at few desired observation points after certain (fairly long) 
interval periodically. 

The LOD2 algorithm above, with its facilitation 
for parallel/reduced/infrequent output processing, 
has been found to be more efficient than those of ADI and SS2 methods 
in their original implementations \cite{LOD3D}. 
To improve the efficiency further, 
we perform simplifications like previous subsections 
to arrive at matrix-operator-free right-hand sides 
for the main iterations as 
\begin{subequations} 
\label{fastLOD2}
\begin{gather}
   \Big( \half \IMat - \fracdt{4} \AMat \Big) \vCol^{n+\tquar}
   = \uCol^{n+\quar}
   \label{fastLOD2a}
\\
   \uCol^{n+\tquar} 
   = \vCol^{n+\tquar} - \uCol^{n+\quar} 
   \label{fastLOD2b}
\\      
   \Big( \half \IMat - \fracdt{4} \BMat \Big) \vCol^{n+\oquar} 
   = \uCol^{n+\tquar}
   \label{fastLOD2c}
\\
   \uCol^{n+\oquar} 
   = \vCol^{n+\oquar} - \uCol^{n+\tquar} 
   \label{fastLOD2d}
.\end{gather}
\end{subequations}
If desired, similar simplifications may be performed for the input processing
\begin{subequations}
\label{fastLOD2i}
\begin{gather}
   \Big( \half \IMat - \fracdt{8} \BMat \Big) \vCol^{\quar} 
   = \uCol^{0}
\\
   \uCol^{\quar} 
   = \vCol^{\quar} - \uCol^{0} 
\end{gather}
\end{subequations}
and output processing 
\begin{subequations}
\label{fastLOD2o}
\begin{gather}
   \Big( \half \IMat + \fracdt{8} \BMat \Big) \vCol^{n+1} 
   = \uCol^{n+\oquar} 
\\
   \uCol^{n+1} 
   = \vCol^{n+1} - \uCol^{n+\oquar} 
.\end{gather}
\end{subequations}
Equations \eqref{fastLOD2a}-\eqref{fastLOD2d} constitute 
the most efficient LOD2 scheme which is comparable to that of 
\eqref{fasterADIa}-\eqref{fasterADId}.

\subsection{Fundamental Implications and Significance}
\label{validext}

It turns out that the new algorithm implementations 
presented above represent 
some very fundamental ones for implicit finite-difference methods. 
As is evident from \eqref{fasterADI} and 
\eqref{fastLOD} or \eqref{fastLOD2}, 
these schemes merely involve similar updating structures, namely 
implicit updating with (simplest) matrix-operator-free right-hand sides 
and explicit updating via subtraction of two vectors. 
Their main slight difference might just be 
the sequence of implicit and explicit updatings. 
Moreover, one may convert from one scheme to another 
simply by proper initialization  
and swapping the roles of field and auxiliary variables. 
For instance, treating $\vCol^{n}$ as the field variables 
(uncontaminated to second order) in \eqref{fastLOD} and using 
\begin{gather}
   \uCol^{n}
   = \Big( \half \IMat + \fracdt{4} \BMat \Big) \vCol^{n}
\end{gather}
to initialize its iterations, 
one may obtain the $\vCol^{n+1}$ field solutions from \eqref{fastLODc} 
that correspond to the $\uCol^{n+1}$ field solutions from 
\eqref{fasterADI}-\eqref{fasterADIo}. 
Such interesting but subtle link between LOD and ADI schemes 
has been pointed out earlier \cite{linkLODADI69}. 
Through the fundamental schemes herein,  
the link becomes particularly obvious 
with their similar updating structures 
and with the aid of variables 
that are otherwise not defined in the original schemes, 
e.g. $\vCol$ is not directly seen in \eqref{orgLOD}. 

The implications of the fundamental schemes above  
are actually more far-reaching. 
In particular, we note that many other classical 
and recent non-dissipative splitting schemes, 
such as D'Yakonov scheme, delta formulation, 
Crank-Nicolson direct-splitting method etc.  
may be cast into the same simplest form of \eqref{fasterADI} 
featuring matrix-operator-free right-hand sides. 
Such analyses when extended to these schemes 
with distinctive splitting formulae,  
will give further insights into the significance of 
fundamental schemes. 
Let us illustrate this point using the scheme with splitting formulae 
(cf.~generalized matricization of D'Yakonov or Beam-Warming scheme)
\begin{subequations}
\label{orgDY}
\begin{gather}
   \Big( \IMat - \fracdt{2} \AMat \Big) \uCol^{*} 
   = \Big( \IMat + \fracdt{2} \AMat \Big)
     \Big( \IMat + \fracdt{2} \BMat \Big) \uCol^{n} 
\\
   \Big( \IMat - \fracdt{2} \BMat \Big) \uCol^{n+1}
   = \uCol^{*} 
.\end{gather}
\end{subequations}
Such scheme may be reduced to an efficient one 
with the iterations algorithm read as 
\begin{subequations}
\label{fastDY}
\begin{gather}
   \vCol^{n} 
   = \uCol^{n} - \vCol^{n-\half} 
\\
   \Big( \half \IMat - \fracdt{4} \AMat \Big) \uCol^{n+\half} 
   = \vCol^{n}
\\
   \vCol^{n+\half} 
   = \uCol^{n+\half} - \vCol^{n} 
\\      
   \Big( \half \IMat - \fracdt{4} \BMat \Big) \uCol^{n+1}
   = \vCol^{n+\half}
\end{gather}
\end{subequations}
along with input initialization
\begin{gather}
   \vCol^{-\half} = \Big( \half \IMat - \fracdt{4} \BMat \Big) \uCol^{0}
.\end{gather}
Notice that the form of \eqref{fastDY} is just that of \eqref{fasterADI}, 
while noting also the relation 
\begin{gather}
   \uCol^{*} = 2 \vCol^{n+\half}  
.\end{gather}
Another illustration is the splitting formulae 
(cf.~generalized matricization of 2D Douglas-Gunn scheme or delta formulation) 
\begin{subequations}
\begin{gather}
   \Big( \IMat - \fracdt{2} \AMat \Big) \Delta \uCol^{*} 
   = \Delta t \Big( \AMat + \BMat \Big) \uCol^{n} 
\\
   \Big( \IMat - \fracdt{2} \BMat \Big) \Delta \uCol^{}
   = \Delta \uCol^{*} 
\\
   \uCol^{n+1} = \uCol^{n} + \Delta \uCol^{}
.\end{gather}
\end{subequations}
Again, this can be turned into the previous fundamental scheme, 
specifically \eqref{uTldCol}-\eqref{fastADIi} and with 
\begin{subequations}
\label{fastDGo}
\begin{gather}
   \Delta \uCol^{*} = \uTldCol^{n+\half} - \uTldCol^{n}
\\
   \Delta \uCol^{} = \half \Big( \uTldCol^{n+1} - \uTldCol^{n} \Big) 
.\end{gather}
\end{subequations}

Meanwhile, the fundamental scheme of \eqref{fastLODa}-\eqref{fastLODb} 
forms the basis of simplification for many split-step methods of 
higher order accuracy in \cite{splitstep03}. 
This will lead to their convenient algorithms 
all with matrix-operator-free right-hand sides 
following the same way that underpins 
\eqref{fastLODc}-\eqref{fastLODd}, \eqref{fastSS2a}-\eqref{fastSS2f} 
and \eqref{fastLOD2a}-\eqref{fastLOD2d}, etc. 
Even for the classical Crank-Nicolson scheme 
having the original unfactorized form 
\begin{gather}
   \bigg[ \IMat - \fracdt{2} \Big( \AMat + \BMat \Big) \bigg] \uCol^{n+1} 
   = \bigg[ \IMat + \fracdt{2} \Big( \AMat + \BMat \Big) \bigg] \uCol^{n} 
,\end{gather}
one may based on similar manipulation to devise 
\begin{subequations}
\begin{gather}
   \bigg[ \half \IMat - \fracdt{4} \Big( \AMat + \BMat \Big) \bigg] \uCol^{n+\half} 
   = \uCol^{n}
\\
   \uCol^{n+1} 
   = \uCol^{n+\half} - \uCol^{n} 
.\end{gather}
\end{subequations}
This obviously constitutes the simpler algorithm 
than that of the original scheme.  

It should be mentioned that the extent of simplification and 
actual improvement of computation efficiency 
for various fundamental schemes 
depend much on the particular matrix operators 
and algorithm implementation details (to be described next).   
Although the scope of this paper is mainly about 
the unconditionally stable implicit FDTD methods in electromagnetics, 
it is evident that all fundamental schemes discussed 
can be extended readily to many other finite-difference schemes 
in computational physics, even conditionally stable  
or parabolic/elliptic ones, etc. 
Some numerical methods that bear resemblance 
to the schemes above such as implicit integrations 
or iterative solutions of certain control problems 
may be made simpler and more efficient in the similar manner.

\section{Efficient Implementations}
\label{implementation}

In this section, we describe the detailed algorithms 
for new efficient implementations of 
the unconditionally stable implicit FDTD methods  
based on the fundamental schemes. 
For simplicity, consider the lossless isotropic medium 
with permittivity $\epsilon$ and permeability $\mu$. 
The splitting matrix operators of Maxwell's equations 
may be selected as 
\begin{gather}
   \AMat = 
   \begin{bmatrix}
      0 & 0 & 0 & 0 & 0 & \frac{1}{\epsilon} \delTld_y
   \\
      0 & 0 & 0 & \frac{1}{\epsilon} \delTld_z & 0 & 0 
   \\
      0 & 0 & 0 & 0 & \frac{1}{\epsilon} \delTld_x & 0 
   \\
      0 & \frac{1}{\mu} \delTld_z & 0 & 0 & 0 & 0   
   \\
      0 & 0 & \frac{1}{\mu} \delTld_x & 0 & 0 & 0   
   \\
      \frac{1}{\mu} \delTld_y & 0 & 0 & 0 & 0 & 0   
   \end{bmatrix}
   \label{AMat}
\\
   \BMat = 
   \begin{bmatrix}
      0 & 0 & 0 & 0 & \frac{-1}{\epsilon} \delTld_z & 0 
   \\
      0 & 0 & 0 & 0 & 0 & \frac{-1}{\epsilon} \delTld_x 
   \\
      0 & 0 & 0 & \frac{-1}{\epsilon} \delTld_y & 0 & 0 
   \\
      0 & 0 & \frac{-1}{\mu} \delTld_y & 0 & 0 & 0   
   \\
      \frac{-1}{\mu} \delTld_z & 0 & 0 & 0 & 0 & 0   
   \\
      0 & \frac{-1}{\mu} \delTld_x & 0 & 0 & 0 & 0   
   \end{bmatrix}
\end{gather}
where $\delTld_x$, $\delTld_y$, $\delTld_z$ 
are the spatial difference operators for the first derivatives 
along $x$, $y$, $z$ directions respectively. 
The vector components of field and auxiliary variables 
can be written as 
\begin{gather}
   \uCol = 
   \begin{bmatrix}
      E_x \\ E_y \\ E_z \\ H_x \\ H_y \\ H_z
   \end{bmatrix}
, \quad 
   \uTldCol = 
   \begin{bmatrix}
      \ETld_x \\ \ETld_y \\ \ETld_z \\ \HTld_x \\ \HTld_y \\ \HTld_z
   \end{bmatrix}
, \quad 
   \vCol = 
   \begin{bmatrix}
      e_x \\ e_y \\ e_z \\ h_x \\ h_y \\ h_z
   \end{bmatrix}
.\end{gather}
For convenience, we also introduce the notations 
\begin{gather}
   b = \fracdt{2 \epsilon}
, \quad
   d = \fracdt{2 \mu}
   \label{bd}
.\end{gather}

\subsection{ADI}
\label{implementADI}

Let us illustrate the implementation of fundamental ADI scheme 
\eqref{fasterADI}. 
While equation~\eqref{fasterADIa} can be implemented directly as 
\begin{subequations}
\label{fasterADIaa}
\begin{gather}
   e_\xi^{n} 
   = \ETld_\xi^{n} - e_\xi^{n-\half}
   , \quad \xi=x,y,z
\\
   h_\xi^{n} 
   = \HTld_\xi^{n} - h_\xi^{n-\half}
   , \quad \xi=x,y,z
   \label{fasterADIaah}
,\end{gather}
\end{subequations}
equation~\eqref{fasterADIb} may find more convenience 
through implicit updating
\begin{subequations}
\begin{gather}
   \half \ETld_x^{n+\half} - \frac{b d}{2} \delTld_y^2 \ETld_x^{n+\half}
   = e_x^{n} + b \delTld_y h_z^{n} 
\\
   \half \ETld_y^{n+\half} - \frac{b d}{2} \delTld_z^2 \ETld_y^{n+\half}
   = e_y^{n} + b \delTld_z h_x^{n} 
\\
   \half \ETld_z^{n+\half} - \frac{b d}{2} \delTld_x^2 \ETld_z^{n+\half}
   = e_z^{n} + b \delTld_x h_y^{n} 
,\end{gather}
\end{subequations}
followed by explicit updating 
\begin{subequations}
\label{fasterADIbbb}
\begin{gather}
   \HTld_x^{n+\half}
   = 2 h_x^{n} + d \delTld_z \ETld_y^{n+\half}
\\
   \HTld_y^{n+\half}
   = 2 h_y^{n} + d \delTld_x \ETld_z^{n+\half}
\\
   \HTld_z^{n+\half}
   = 2 h_z^{n} + d \delTld_y \ETld_x^{n+\half}
.\end{gather}
\end{subequations}
Similar arguments apply to equations~\eqref{fasterADIc} and 
\eqref{fasterADId}, where the former is implemented directly as 
\begin{subequations}
\label{fasterADIcc}
\begin{gather}
   e_\xi^{n+\half} 
   = \ETld_\xi^{n+\half} - e_\xi^{n}
   , \quad \xi=x,y,z
\\
   h_\xi^{n+\half} 
   = \HTld_\xi^{n+\half} - h_\xi^{n}
   , \quad \xi=x,y,z
   \label{fasterADIcch}
,\end{gather}
\end{subequations}
while the latter comprises implicit updating
\begin{subequations}
\label{fasterADIdd}
\begin{gather}
   \half \ETld_x^{n+1} - \frac{b d}{2} \delTld_z^2 \ETld_x^{n+1}
   = e_x^{n+\half} - b \delTld_z h_y^{n+\half} 
\\
   \half \ETld_y^{n+1} - \frac{b d}{2} \delTld_x^2 \ETld_y^{n+1}
   = e_y^{n+\half} - b \delTld_x h_z^{n+\half} 
\\
   \half \ETld_z^{n+1} - \frac{b d}{2} \delTld_y^2 \ETld_z^{n+1}
   = e_z^{n+\half} - b \delTld_y h_x^{n+\half} 
,\end{gather}
\end{subequations}
and explicit updating
\begin{subequations}
\label{fasterADIddd}
\begin{gather}
   \HTld_x^{n+1}
   = 2 h_x^{n+\half} - d \delTld_y \ETld_z^{n+1}
\\
   \HTld_y^{n+1}
   = 2 h_y^{n+\half} - d \delTld_z \ETld_x^{n+1}
\\
   \HTld_z^{n+1}
   = 2 h_z^{n+\half} - d \delTld_x \ETld_y^{n+1}
.\end{gather}
\end{subequations}

When the difference operators above represent specifically 
the second-order central-differencing operators on Yee cells, 
equations~\eqref{fasterADIaa}-\eqref{fasterADIddd} 
correspond to the efficient algorithm delineated in \cite{ADIefficient}. 
In particular, the solution $\ETld_\xi^{n+1}$ from \eqref{fasterADIdd} 
coincides with $E_\xi^*|^{n+1}$ from \cite[eqn.~(11)]{ADIefficient}, 
while the solution $\HTld_\xi^{n+1}$ from \eqref{fasterADIddd} 
is half of $H_\xi^*|^{n+1}$ from \cite[eqn.~(12)]{ADIefficient}. 
The reason for such magnetic field relation is 
that certain variables and coefficients in 
\eqref{fasterADIaa}-\eqref{fasterADIddd} have been re-scaled  
in the corresponding auxiliary and explicit updating equations of 
\cite{ADIefficient}. 
This helps save the scaling of final magnetic field solution 
that is otherwise necessary, cf.~\eqref{fasterADIo}. 
Moreover, additional savings of operations can be achieved 
by combining \eqref{fasterADIbbb} and \eqref{fasterADIcch} as 
\begin{subequations}
\label{fasterADIbbh}
\begin{gather}
   h_x^{n+\half}
   = h_x^{n} + d \delTld_z \ETld_y^{n+\half}
\\
   h_y^{n+\half}
   = h_y^{n} + d \delTld_x \ETld_z^{n+\half}
\\
   h_z^{n+\half}
   = h_z^{n} + d \delTld_y \ETld_x^{n+\half}
.\end{gather}
\end{subequations}
When the final magnetic field data is not needed frequently, 
\eqref{fasterADIddd} and \eqref{fasterADIaah} (at $n+1$ time step) 
may also be combined for higher efficiency: 
\begin{subequations}
\label{fasterADIddh}
\begin{gather}
   h_x^{n+1}
   = h_x^{n+\half} - d \delTld_y \ETld_z^{n+1}
\\
   h_y^{n+1}
   = h_y^{n+\half} - d \delTld_z \ETld_x^{n+1}
\\
   h_z^{n+1}
   = h_z^{n+\half} - d \delTld_x \ETld_y^{n+1}
.\end{gather}
\end{subequations}
Other implementation details may be referred to \cite{ADIefficient} 
including for-looping, tridiagonal system solving, memory reuse etc. 
In actuality, there are many possibilities 
as far as detailed implementations are concerned, 
such as different definition and scaling of variables, 
different implicit and explicit updatings, 
different difference operators, 
e.g. higher order, 
parameter optimized, 
compact 
etc. 
\cite{FuADIhigh05}-\cite{optimized03}.
The updating equations presented above still leave much 
room and generality for exploring into these possibilities.

\subsection{LOD2}
\label{implementLOD2}

For subsequent discussions and comparisons, 
we also describe the algorithm implementation of 
fundamental LOD2 scheme \eqref{fastLOD2} 
in more detail and yet general. 
Adopting the same notations of \eqref{AMat}-\eqref{bd}, 
equation~\eqref{fastLOD2a} may be split into implicit updating 
\begin{subequations}
\label{fastLOD2aa}
\begin{gather}
   \half e_x^{n+\tquar} - \frac{b d}{2} \delTld_y^2 e_x^{n+\tquar}
   = E_x^{n+\quar} + b \delTld_y H_z^{n+\quar} 
   \label{fastLOD2.1}
\\
   \half e_y^{n+\tquar} - \frac{b d}{2} \delTld_z^2 e_y^{n+\tquar}
   = E_y^{n+\quar} + b \delTld_z H_x^{n+\quar} 
   \label{fastLOD2.2}
\\
   \half e_z^{n+\tquar} - \frac{b d}{2} \delTld_x^2 e_z^{n+\tquar}
   = E_z^{n+\quar} + b \delTld_x H_y^{n+\quar} 
   \label{fastLOD2.3}
\end{gather}
\end{subequations}
and explicit updating 
\begin{subequations}
\label{fastLOD2aaa}
\begin{gather}
   h_x^{n+\tquar}
   = 2 H_x^{n+\quar} + d \delTld_z e_y^{n+\tquar}
   \label{fastLOD2.4}
\\
   h_y^{n+\tquar}
   = 2 H_y^{n+\quar} + d \delTld_x e_z^{n+\tquar}
   \label{fastLOD2.5}
\\
   h_z^{n+\tquar}
   = 2 H_z^{n+\quar} + d \delTld_y e_x^{n+\tquar}
   \label{fastLOD2.6}
,\end{gather}
\end{subequations}
while \eqref{fastLOD2b} is directly 
\begin{subequations}
\begin{gather}
   E_\xi^{n+\tquar} 
   = e_\xi^{n+\tquar} - E_\xi^{n+\quar}
   , \quad \xi=x,y,z
   \label{fastLOD2.bbE}
\\
   H_\xi^{n+\tquar} 
   = h_\xi^{n+\tquar} - H_\xi^{n+\quar}
   , \quad \xi=x,y,z
   \label{fastLOD2.bbH}
.\end{gather}
\end{subequations}
Likewise, equation~\eqref{fastLOD2c} consists of implicit updating
\begin{subequations}
\begin{gather}
   \half e_x^{n+\oquar} - \frac{b d}{2} \delTld_z^2 e_x^{n+\oquar}
   = E_x^{n+\tquar} - b \delTld_z H_y^{n+\tquar} 
   \label{fastLOD2.13}
\\
   \half e_y^{n+\oquar} - \frac{b d}{2} \delTld_x^2 e_y^{n+\oquar}
   = E_y^{n+\tquar} - b \delTld_x H_z^{n+\tquar} 
   \label{fastLOD2.14}
\\
   \half e_z^{n+\oquar} - \frac{b d}{2} \delTld_y^2 e_z^{n+\oquar}
   = E_z^{n+\tquar} - b \delTld_y H_x^{n+\tquar} 
   \label{fastLOD2.15}
\end{gather}
\end{subequations}
and explicit updating 
\begin{subequations}
\label{fastLOD2ccc}
\begin{gather}
   h_x^{n+\oquar}
   = 2 H_x^{n+\tquar} - d \delTld_y e_z^{n+\oquar}
   \label{fastLOD2.16}
\\
   h_y^{n+\oquar}
   = 2 H_y^{n+\tquar} - d \delTld_z e_x^{n+\oquar}
   \label{fastLOD2.17}
\\
   h_z^{n+\oquar}
   = 2 H_z^{n+\tquar} - d \delTld_x e_y^{n+\oquar}
   \label{fastLOD2.18}
,\end{gather}
\end{subequations}
whereas \eqref{fastLOD2d} is simply
\begin{subequations}
\label{fastLOD2dd}
\begin{gather}
   E_\xi^{n+\oquar} 
   = e_\xi^{n+\oquar} - E_\xi^{n+\tquar}
   , \quad \xi=x,y,z
   \label{fastLOD2.ddE}
\\
   H_\xi^{n+\oquar} 
   = h_\xi^{n+\oquar} - H_\xi^{n+\tquar}
   , \quad \xi=x,y,z
   \label{fastLOD2.ddH}
.\end{gather}
\end{subequations}

\begin{table*}[t]
\caption{Comparisons of Unconditionally Stable FDTD Methods with Second-order Spatial Central Difference}
\centering
\begin{tabular}{|c|c|c|c|c|c|c|c|c|c|}
\hline
\multicolumn{2}{|c|}{Scheme} 
& \multicolumn{2}{|c|}{ADI} 
& \multicolumn{2}{|c|}{LOD1/SS1}
& \multicolumn{2}{|c|}{SS2}
& \multicolumn{2}{|c|}{LOD2}
\\
\hline
\multicolumn{2}{|c|}{Algorithm} 
& Original & New 
& Original & New 
& Original & New 
& Original & New 
\\
\hline
\multirow{3}{10ex}{Equations} 
& Main Iterations 
& \eqref{clsADI} & \eqref{fasterADI}
& \eqref{orgLOD} & \eqref{fastLOD} 
& \eqref{orgSS2} & \eqref{fastSS2} 
& \eqref{orgLOD2} & \eqref{fastLOD2} 
\\
\cline{2-10}
& Input & -- & \eqref{uTldCol}, \eqref{fastADIi} & -- & -- & -- & -- & \eqref{orgLOD2i} & \eqref{fastLOD2i} 
\\
\cline{2-10}
& Output & -- & \eqref{fasterADIo} & -- & -- & -- & -- & \eqref{orgLOD2o} & \eqref{fastLOD2o} 
\\
\hline
\multirow{4}{10ex}{Memory} 
& Field Array 1
& \parbox{8ex}
{$\uCol^{n+1} \\ \uCol^{n}$}
& \parbox{8ex}
{$\uTldCol^{n+1} \\ \vCol^{n} \\ \uTldCol^{n}$}
& \parbox{8ex}
{$\uCol^{n+1} \\ \uCol^{n}$}
& \parbox{8ex}
{$\uCol^{n+1} \\ \vCol^{n+1} \\ \uCol^{n}$}
& \parbox{8ex}
{$\uCol^{n+1} \\ \uCol^{n+\quar}$}
& \parbox{8ex}
{$\uCol^{n+1} \\ \vCol^{n+1} \\ \uCol^{n+\quar} \\ \vCol^{n+\quar}$}
& \parbox{8ex}
{$\uCol^{n+\oquar} \\ \uCol^{n+\quar}$}
& \parbox{8ex}
{$\uCol^{n+\oquar} \\ \vCol^{n+\oquar} \\ \uCol^{n+\quar}$}
\\
\cline{2-10}
& Field Array 2
& \parbox{8ex}
{$\uCol^{n+\half}$} 
& \parbox{8ex}
{$\vCol^{n+\half} \\ \uTldCol^{n+\half} \\ \vCol^{n-\half}$} 
& \parbox{8ex}
{$\uCol^{n+\half}$} 
& \parbox{8ex}
{$\uCol^{n+\half} \\ \vCol^{n+\half}$} 
& \parbox{8ex}
{$\uCol^{n+\tquar} \\ \uCol^{n}$} 
& \parbox{8ex}
{$\uCol^{n+\tquar} \\ \vCol^{n+\tquar} \\ \uCol^{n}$} 
& \parbox{8ex}
{$\uCol^{n+\tquar}$} 
& \parbox{8ex}
{$\uCol^{n+\tquar} \\ \vCol^{n+\tquar}$} 
\\
\hline
\multirow{2}{10ex}{Implicit} 
& M/D & 18 & 6 & 18 & 6 & 27 & 9 & 18 & 6
\\
\cline{2-10}
& A/S & 48 & 18 & 24 & 18 & 36 & 27 & 24 & 18 
\\
\hline
\multirow{2}{10ex}{Explicit} 
& M/D & 12 & 6 & 6 & 6 & 9 & 9 & 6 & 6
\\
\cline{2-10}
& A/S & 24 & 12 & 24 & 12 & 36 & 18 & 24 & 12  
\\
\hline
\multirow{3}{10ex}{Total} 
& M/D & 30 & 12 & 24 & 12 & 36 & 18 & 24 & 12  
\\
\cline{2-10}
& A/S & 72 & 30 & 48 & 30 & 72 & 45 & 48 & 30
\\
\cline{2-10}
& M/D+A/S & 102 & 42 & 72 & 42 & 108 & 63 & 72 & 42 
\\
\hline
\multirow{2}{10ex}{Efficiency gain} 
& RHS & 1 & 2.43 & 1.42 & 2.43 & 0.94 & 1.62 & 1.42 & 2.43
\\
\cline{2-10}
& Overall & 1 & 1.83 & 1.29 & 1.83 & 0.86 & 1.22 & 1.29 & 1.83  
\\
\hline
\multicolumn{2}{|c|}{For-Loops} 
& \multicolumn{2}{|c|}{12} 
& \multicolumn{2}{|c|}{12} 
& \multicolumn{2}{|c|}{18} 
& \multicolumn{2}{|c|}{12} 
\\
\hline
\multicolumn{2}{|c|}{Temporal Accuracy} 
& \multicolumn{2}{|c|}{Second-Order} 
& \multicolumn{2}{|c|}{First-Order} 
& \multicolumn{2}{|c|}{Second-Order} 
& \multicolumn{2}{|c|}{Second-Order} 
\\
\hline
\end{tabular}
\label{flops}
\end{table*}

As before, there are many choices for the difference operators 
in equations~\eqref{fastLOD2aa}-\eqref{fastLOD2dd}. 
For being more specific in the next section, we select 
the second-order central-differencing operators on Yee cells. 
To minimize the number of for-loops in the previous algorithm, 
some of their updating equations may be incorporated 
in the same loops. 
Moreover, when the auxiliary magnetic variables $h_\xi$ 
are not to be output, the explicit updating equations 
\eqref{fastLOD2aaa} and \eqref{fastLOD2.bbH} can be combined to read 
\begin{subequations}
\label{fastestLOD2ex1}
\begin{gather}
   H_x^{n+\tquar}
   = H_x^{n+\quar} + d \delTld_z e_y^{n+\tquar}
\\
   H_y^{n+\tquar}
   = H_y^{n+\quar} + d \delTld_x e_z^{n+\tquar}
\\
   H_z^{n+\tquar}
   = H_z^{n+\quar} + d \delTld_y e_x^{n+\tquar}
,\end{gather}
\end{subequations}
while \eqref{fastLOD2ccc} and \eqref{fastLOD2.ddH} become 
\begin{subequations}
\label{fastestLOD2ex2}
\begin{gather}
   H_x^{n+\oquar}
   = H_x^{n+\tquar} - d \delTld_y e_z^{n+\oquar}
\\
   H_y^{n+\oquar}
   = H_y^{n+\tquar} - d \delTld_z e_x^{n+\oquar}
\\
   H_z^{n+\oquar}
   = H_z^{n+\tquar} - d \delTld_x e_y^{n+\oquar}
.\end{gather}
\end{subequations}
Meanwhile, to minimize the memory storage requirement, 
some of the variables may occupy the same memory spaces. 
For instance, one can choose to reuse the spaces of 
those variables grouped within curly braces as ($\xi=x,y,z$): 
\begin{gather*}
   \{E_\xi^{n+\quar},e_\xi^{n+\oquar},E_\xi^{n+\oquar}\}; 
\quad
   \{e_\xi^{n+\tquar},E_\xi^{n+\tquar}\};
\\
   \{H_\xi^{n+\quar},h_\xi^{n+\oquar},H_\xi^{n+\oquar}\};
\quad
   \{h_\xi^{n+\tquar},H_\xi^{n+\tquar}\}
.\end{gather*}
In a computer program, these variables could just 
share the same names to be assigned with new values successively. 
Note that in conjunction with 
\eqref{fastestLOD2ex1}-\eqref{fastestLOD2ex2},
one may absorb $h_\xi$ altogether and reuse the $H_\xi$ spaces. 

Following the similar arguments as above, 
the detailed algorithms for the fundamental LOD1/SS1 and SS2 schemes 
may be written and coded with reference to 
\eqref{fastLOD} and \eqref{fastSS2} accordingly.

\section{Discussions and Comparisons}
\label{comparison}

Having systematically addressed various implicit schemes 
and their new algorithm implementations, we are now ready to make 
a comparative study. 
Of particular interest are the computation efficiency gains 
as compared to the original schemes, or more specifically 
over the prevailing ADI-FDTD implementation 
\cite{ADI00, Namiki00}. 
Furthermore, the fundamental nature of new algorithm implementations,  
as has been pointed out in \secref{validext}, 
will become evident once again 
through the comparisons and discussions. 

In what follows, we acquire the floating-point operations 
(flops) count for the main iterations of 
ADI, LOD1/SS1, SS2 and LOD2 schemes. 
Based on the right-hand sides (RHS) of their respective 
updating equations using second-order central-differencing operators, 
cf.~Sections~\ref{formulation} and \ref{implementation},
the number of multiplications/divisions (M/D) and 
additions/subtractions (A/S) required 
for one complete time step are determined. 
Each complete step comprises two/three procedures 
whose implicit and explicit updating codes 
are assumed to have been arranged 
in order to achieve minimum number of for-loops. 
It is also assumed that all multiplicative factors 
have been precomputed and stored, 
while the number of electric and magnetic field components 
in all directions have been taken to be the same. 
Our results are summarized in Table~\ref{flops}, 
which lists the flops count for both original  
and new algorithm implementations of each scheme. 
For convenience, the pertaining equations involved 
in the main iterations are also indicated, 
along with those that may be needed for 
input initialization and output processing. 
When there is no input or output equation shown, 
it means that no special input or output treatment 
is required for the particular implementation. 
For the main iterations variables, 
their memory storage requirements are also specified 
in Table~\ref{flops}. 
By properly reusing the spaces of those variables 
grouped within the same table entity, 
only two field arrays need to be allocated for all. 
Note that judicious reuse of memory space is important 
so that the iterations do not often invoke virtual memory 
(e.g. via hard disk). 
Otherwise, they may lead to 
slower computations than those iterations 
with more flops but do not need virtual memory. 

From Table~\ref{flops}, it is clear that 
the total flops count (M/D+A/S) for all implicit schemes 
has been reduced considerably using the new implementations 
as compared to the original counterparts. 
Most notably, the conventional ADI-FDTD method based on 
\eqref{clsADI} takes 102 flops \cite{ADI00, Namiki00}, 
whereas the new one based on \eqref{fasterADI} 
merely takes 42 flops, cf.~\secref{implementADI} 
(assuming no need for frequent output of magnetic field data). 
Notice that the LOD1/SS1 and LOD2 schemes, which take 
72 flops in their original implementations \cite{LOD3D}, 
also take the same 42 flops using the new algorithms, 
cf.~\secref{implementLOD2}. 
Moreover, the number and type of flops for respective 
implicit and explicit updatings are seen to be 
identical for all these new implicit schemes. 
This is simply a manifestation of their aforementioned 
fundamental updating structures, 
which can be found to take 21 flops per updating procedure.  
Applying such argument for the SS2 scheme 
with three updating procedures, it is easy 
to deduce that its new algorithm will take 63 flops, 
in contrast to the 108 flops required 
when using its original implementation. 

Based on the flops count reduction 
for the right-hand sides of updating equations, 
Table~\ref{flops} has listed the efficiency gains 
over the conventional ADI-FDTD method. 
The new fundamental ADI, LOD1/SS1 and LOD2 schemes 
all feature the same efficiency gain of 2.43 
with their same reduced total flops count. 
Since there is also the cost for solving tridiagonal systems, 
we estimate about $5N$ flops for a system of order $N$ 
using precomputed bidiagonally factorized elements.  
Taking this cost into account, all these new implementations 
still achieve an overall efficiency gain of 1.83  
in flops count reduction over the conventional ADI-FDTD. 
It is interesting to find that 
even for the SS2 scheme with three updating procedures, 
the overall efficiency may still be higher (gain $\sim 1.22$) 
by using the new algorithm implementation. 
Note that although we have carefully taken into account 
various flops count, 
there still exist other factors that may affect 
the actual computation efficiency. 
These factors are often dependent on the particular 
computer platform, hardware configuration, operating system, 
compiler software, program code arrangement, etc. 
\cite{SunTrueman05}.

Besides the flops count, 
the for-loops count is also given in Table~\ref{flops}. 
Each for-loop is to perform the whole sweep of 
$i$, $j$, $k$ indices along $x$, $y$, $z$ directions 
for one field component. 
For both original and new implementations of 
ADI, LOD1/SS1 and LOD2 schemes, 
there are 12 for-loops for all updating equations 
in both procedures. 
For the SS2 scheme, there are 18 for-loops altogether 
corresponding to its three updating procedures. 
Meanwhile, the order of temporal accuracy is also summarized 
for various implicit schemes in Table~\ref{flops}. 
As mentioned in \secref{validext}, even though 
the LOD1/SS1 scheme is only first-order accurate in time, 
it may be converted to other schemes of second-order accuracy 
by exploiting their similar fundamental updating structures. 

Although not included in Table~\ref{flops}, 
many other implicit schemes discussed 
in \secref{validext} 
can be found to take the same reduced flops count (42)
and the same for-loops count (12) when using 
the new implementations with \eqref{AMat}-\eqref{bd}.  
This is anticipated because they all simply converge 
to the same fundamental schemes 
with matrix-operator-free right-hand sides 
just like those addressed in the table, 
even though their original splitting formulae are distinctive, 
cf.~\eqref{orgDY}-\eqref{fastDGo}. 
Here again the significance of these fundamental schemes accrues 
when one is equipped with a general methodical approach 
for improving the computation efficiency of 
their original counterparts. 
Moreover, their updating structures that are of fundamental nature 
make them well-suited to serve as the basis and benchmark 
for construction and development of future implicit schemes, 
e.g. with new $\AMat$ and $\BMat$ 
or different time-stepping procedures, etc.

\section{Conclusion} 
This paper has presented the generalized formulations of 
fundamental schemes 
for efficient unconditionally stable implicit FDTD methods. 
The formulations have been presented in terms of 
generalized matrix operator equations 
pertaining to classical splitting formulae 
of ADI, LOD and split-step schemes.
To provide further insights into the implications 
and significance of fundamental schemes, 
the analyses have also been extended to many other schemes 
with distinctive splitting formulae. 
It has been noted that all the fundamental schemes 
feature similar fundamental updating structures 
that are in simplest forms 
with most efficient right-hand sides. 
Detailed algorithms have been described 
for new efficient implementations of 
the unconditionally stable implicit FDTD methods 
based on the fundamental schemes. 
A comparative study of various implicit schemes 
in their original and new implementations has been carried out,  
which includes comparisons of their computation costs 
and efficiency gains.  

The fundamental schemes presented in this paper will be 
of much usefulness and significance not only in electromagnetics, 
but also in many other areas that may be adopting 
various classical implicit schemes.  
With their simplest forms featuring 
the most efficient right-hand sides, 
these fundamental schemes will lead to coding simplification 
and efficiency improvement in algorithm implementations. 
Furthermore, their fundamental updating structures 
invite further investigations into their properties 
and subsequent extensions for their applications. 
They will also serve aptly as the basis and benchmark 
for construction and development of future implicit schemes.



\end{document}